\documentstyle[floats,prl,aps]{revtex}
\begin{document}
\renewcommand{\thefootnote}{\fnsymbol{footnote}}
\draft
                     


\def\a{\alpha}
\def\b{\beta}
\def\d{\delta}
\def\e{\epsilon}
\def\g{\gamma}
\def\k{\kappa}
\def\l{\lambda}
\def\o{\omega}
\def\t{\theta}
\def\s{\sigma}
\def\D{\Delta}
\def\L{\Lambda}
\def\Uq{U_{q}(\widehat{sl}(2|1))}


\def\beq{\begin{equation}}
\def\eeq{\end{equation}}
\def\bea{\begin{eqnarray}}
\def\eea{\end{eqnarray}}
\def\ba{\begin{array}}
\def\ea{\end{array}}
\def\no{\nonumber}
\def\le{\langle}
\def\re{\rangle}
\def\lt{\left}
\def\rt{\right}

\newtheorem{Theorem}{Theorem}
\newtheorem{Definition}{Definition}
\newtheorem{Proposition}{Proposition}
\newtheorem{Lemma}{Lemma}
\newtheorem{Corollary}{Corollary}
\newcommand{\proof}[1]{{\bf Proof. }
        #1\begin{flushright}$\Box$\end{flushright}}

\newcommand{\sect}[1]{\setcounter{equation}{0}\section{#1}}
\renewcommand{\theequation}{\thesection.\arabic{equation}}

\title{\large\bf
Highest weight representations of $\Uq$ and
correlation functions of the $q$-deformed supersymmetric $t$-$J$ model}
\author{\large Wen-Li Yang$^{~1,2}$ and Yao-Zhong Zhang$^{~2}$}

\address{$^{1}~$ Institute of Modern Physics,
Northwest University, Xian  710069,China\\
$^{2}~$ Department of Mathematics, University of Queensland,
Brisbane, Qld 4072, Australia}

\maketitle
\vspace{10pt}

\begin{abstract}
We re-examine the level-one irreducible highest weight representations
of the quantum affine superalgebra $U_q(\widehat{sl}(2|1))$ and derive
the characters and supercharacters associated with these
representations. We calculate the exchange relations of the
vertex operators and find that these vertex operators satisfy the graded
Faddeev-Zamolodchikov algebra. We give an integral expression of
the correlation functions of the $q$-deformed supersymmetric $t$-$J$ model
and derive the difference equations which they satisfy.
\end{abstract}

\section{Introduction} 

The algebraic analysis method based on 
infinite dimensional non-abelian symmetries such as Virasoro and
Kac-Moody algebra symmetries
has proved eminently successful in both formulating and solving
low-dimensional systems on the critical points (see e.g. \cite{Itz88}
and references therein). The
key elements in this approach are the highest weight 
representation theory and vertex operators which are intertwiners
between two irreducible highest weight representations.
One advantage over other abelian
symmetry methods such as the Bethe ansatz and QISM is that it enables
one to compute correlation functions and form factors exactly in the
form of integral representations. 

Following the success of this approach, one wonders if a similar
program can be carried out for massive intergable systems, i.e.
integrable systems away from the critical points. After the discovery
of quantum groups and quantum affine algebras, one has every reason
to believe that the goal is achievable
because these quantum algebras are exactly the non-abelian symmetries
underlying the off-critical integrable models. One breakthrough is due
to Frenkel and Reshetikhin \cite{Fre92} who introduced the $q$-deformed
vertex operators associated with the quantum affine algebras. They
showed that the correlation functions satisfy a set of difference
equations, the so-called $q$-KZ equations. Using the $q$-vertex
operators and the level one free boson realization of quantum affine
algebras of Frenkel and Jing \cite{Fre88}, the Kyoto group \cite{DFJMN}
(for a nice review, see \cite{JM}) succeeded
in diagonalizing the XXZ spin chain directly in the thermodynamic limit.
The method is very similar to that used in the critical cases. The
idea is to work directly on an infinite lattice and use
the full quantum affine algebra symmetry of the model. This way, 
the highest weight representation
theory and $q$-vertex operators of the non-abelian quantum symmetry 
enter the game in a similar way as in the critical cases. Soon after this
ground-breaking work of the Kyoto group, several generalizations have
been considered (see e.g. \cite{K} for the vertex models
with $U_q(\widehat{sl}(n))$-symmetry and \cite{LP,AJ} for the face type
statistical mechanics models). 

Like in the critical cases,  off-critical integrable models with
(quantum) superalgebra symmetries have occupied an important place. So it
is natural to generalize the above program further to the supersymmetric
case. Quantum (affine) superalgebras have much more
complicated structures and representation theory than their non-super
counterparts \cite{KW}. In general, 
the knowledge about the representation theory
of quantum affine superalgebras is still very much imcomplete. 
For the type I quantum affine superalgebra $U_q(\widehat{sl}
(m|n))$ ($m\neq n$), level-one representations and $q$-vertex operators
have been
studied in \cite{KS} (see \cite{Awa97} for the level-$k$
free boson realization of $\Uq$).  In particular, 
the level-one irreducible highest weight representations of $\Uq$
have been investigated in details and the character formulae for these
modules have been conjectured.

In this paper, we re-examine the level-one irreducible highest weight
representations of $\Uq$ and calculate their characters by means of
the BRST resolution. We also derive the super characters associated
with these modules. We find that the conjecture 2.2 proposed in
\cite{KS} needs to be slightly modified and we give a modified
conjecture (conjecture 1 below). In section 4, we compute the 
exchange relations of the $q$-vertex operators and show that these
vertex operators satisfy the graded Faddeev-Zamolodchikov algebra. 
A Miki's construction of $\Uq$ is also given. In section 5, we
consider the application of the irreducible highest weight modules and
$q$-vertex operators to the $q$-deformed supersymmetric $t$-$J$
model on an infinite lattice. Generalizing the Kyoto group's work
\cite{DFJMN}, we give a mathematical definition of the space of
physical states of the supersymmetric $t$-$J$ model and define its
local structure and local operators. We compute the one-point
correlation functions of the local operators and give an integral
expression of the correlation functions. A set of infinite number of
difference equations satisfied by the correlation functions has also
been derived.

\section{Preliminaries}

In this section, we briefly review the bosonization of quantum affine 
superalgebra $U_q(\widehat{sl}(2|1)$ at level one \cite{KS}.

\subsection{ Quantum affine superalgebra $U_q(\widehat{sl}(2|1)$ }

The symmetric Cartan matrix of the affine Lie superalgebra $\widehat{sl}(2|1)$
is 
\begin{eqnarray*}
(a_{ij})=\left(
\begin{array}{ccc}
0&-1&1\\
-1&2&-1\\
1&-1&0
\end{array}\right)
\end{eqnarray*}
where $i,j=0,1,2$. Quantum affine superalgebra 
$U_q(\widehat{sl}(2|1)$ is a $q$-analogue of the universal 
enveloping algebra of $\widehat{sl}(2|1)$ generated by the Chevalley 
generators $\{e_i,f_i,t_i^{\pm 1},d | i=0,1,2\}$, where $d$ is the
usual derivation operator. The $Z_2$-grading
of the generators are
$[e_0]=[f_0]=[e_{2}]=[f_{2}]=1$ and zero otherwise. 
The defining relations are
\begin{eqnarray*}
& & t_it_j =t_jt_i,\ \  t_id=dt_i, \ \ [d,e_i]=\delta_{i,0}e_i,\ \ 
[d,f_i]=-\delta_{i,0}f_i,\\ 
& &t_ie_jt_i^{-1}=q^{a_{ij}}e_j,\ \ t_if_jt_i^{-1}=q^{-a_{ij}}f_j ,\ \  
[e_i,f_j] =\delta_{ij}{ t_i-t_i^{-1} \over q-q^{-1}},
\end{eqnarray*}
plus the extra $q$-Serre type relations \cite{Y} which we omit.
Here and throughout, $[X,Y]_\xi=XY-(-1)^{[X][Y]}\xi YX$ and $[X,Y]=
[X,Y]_1$.
 
$U_q(\widehat{sl}(2|1))$ is a quasi-triangular Hopf superalgebra
endowed with the ${\bf Z}_2$-graded Hopf algebra structure:
\bea
&&\Delta(t_i)=t_i\otimes t_i,~~~~
   \Delta(e_i)=e_i\otimes 1+t_i\otimes e_i,~~~~
    \Delta(f_i)=f_i\otimes t_i^{-1}+1\otimes f_i ,\no\\
&&\epsilon(t_i)=1,~~~~\epsilon(e_i)=\epsilon(f_i)=0,\no\\
&&S(e_i)=-t_i^{-1} e_i,~~~~ S(f_i)=-f_i t_i, ~~~~
   S(t_i^{\pm 1})=t_i^{\mp 1},~~~~ S(d)=-d.
\eea 
Note the antipode $S$ is a graded algebra anti-automorphism. Namely for
homogeneous elemets $a,b \in \Uq$,
$S(ab)=(-1)^{[a][b]}S(b)S(a)$.
The multiplication rule for the tensor product is ${\bf Z}_2$
graded and is defined for homogeneous elements $a,b,a',b' \in \Uq$ by 
$(a\otimes b)(a'\otimes b')=(-1)^{[b][a']}(aa'\otimes bb')$,  
which extends to inhomgeneous elements through linearity.

$\Uq$ can also be realized by 
the Drinfeld generators \cite{Dri88} $\{d,\ \ X^{\pm,i}_{m}$,
$h^{i}_n$, $K^{i,\pm 1},\gamma^{\pm 1/2} | i=1,2,$ $m \in 
{\bf Z}, n \in {\bf Z}_{\ne 0}\}$. 
The ${\bf Z}_2$-grading of the Drinfeld generators 
are $[X^{\pm,2}_m]=1$ ($m\in {\bf Z}$) and zero otherwise. The
relations read \cite{Y,Z}
\begin{eqnarray*}
& &\gamma\ \  {\rm is\ \  central },\ \
[K^i,h^j_m]=0,\ \ [d,K^i]=0,\ \ [d,h^{j}_m]=mh^j_m\\
& &[h^i_m,h^j_n] =\delta_{m+n,0} 
{[a_{ij}m](\gamma^m-\gamma^{-m}) \over m(q-q^{-1})} \\
& &K^iX^{\pm,j}_m =q^{\pm a_{ij}}X^{\pm,j}_m K^i,\ \
[d,X^{\pm,j}_m]=mX^{\pm,j}_m \\ 
& &[h^i_m,X^{\pm,j}_n]=\pm {[a_{ij}m] \over m}
\gamma^{\pm |m|/2}X^{\pm,j}_{n+m},\\ 
& &[X^{+,i}_m,X^{-,j}_n]=\frac{\delta_{i,j}}{ q-q^{-1}}
(\gamma^{(m-n)/2}\psi^{+,j}_{m+n} -\gamma^{-(m-n)/2 }
\psi^{-,j}_{m+n}),\\ 
& &[X^{\pm ,2}_m,X^{\pm ,2}_n]=0 \\
& &[X^{\pm,i}_{m+1},X^{\pm,j}_{n}]_{q^{\pm a_{ij}}} 
+[X^{\pm,j}_{n+1},X^{\pm,i}_{m}]_{q^{\pm a_{ij}}}=0, 
\ \ {\rm for} \ \ a_{ij}\ne 0,
\end{eqnarray*} 
where $[m]=\frac{q^m-q^{-m}}{q-q^{-1}}$.

The Chevalley generators are related to the Drinfeld generators
by the formulae:
\begin{eqnarray}
& &t_i= K_i,\ \ e_i = X^{+,i}_0,\ \ t_0=\gamma(K^2K^1)^{-1},
\ \ f_i = X^{-,i}_0 \ \  
{\rm for}\ \ i=1,2 ,\\  
& &e_0= -[X_0^{-,2},X^{-,1}_1]_{q{-1}}(K^1K^2)^{-1},\ \ 
f_0=K^1K^2[X^{+,1}_{-1},X^{+,2}_0]_q.
\end{eqnarray}

\subsection{Bosonization of $\Uq$ at level one}
Let us introduce the bosonic $q$-oscillators\cite{KS}
$\{a^1_n,a^2_n,b_n,c_n$, $Q_{a^1},Q_{a^2}, Q_{b},Q_{c}$ 
$|n \in {\bf Z} \}$ which satisfy the commutation relations
\begin{eqnarray*}
& &[a^i_m,a^j_n]=\delta_{i,j}\delta_{m+n,0}{[m]^2 \over m},\ \ \  
[a^i_0,Q_{a^j}] = \delta_{i,j},\\
& &[b_m,b_n]=-\delta_{m+n,0}{[m]^2 \over m},\ \ \  
[b_0,Q_{b}] = -1  \\
& &[c_m,c_n]=\delta_{m+n,0}{[m]^2 \over m},\ \ \ \  
[c_0,Q_{c}] = 1. 
\end{eqnarray*}
\noindent The remaining commutators vanish. 
Define the Drinfeld currents or generating functions 
$X^{\pm}_i(z)=\sum_{m\in {\bf Z}}X^{\pm,i}_mz^{-m-1} $,
and introduce $h^i_0$ by setting $K^i=q^{h^i_0}$.
Let $Q_{h^1}=Q_{a^1}-Q_{a^{2}}$, 
$Q_{h^{2}}=Q_{a^{2}}+Q_{b}$. Define 
\beq
h_i(z;\k)=-\sum_{n \ne 0}{h^i_n \over 
  [n]}q^{-\k|n|}z^{-n} +Q_{h^i}+h^i_0\ln z.
\eeq
Other bosonic fields $b(z;\k)$ and $c(z;\k)$ are defined similarly.
Also introduce the $q$-differential operator defined by 
 $\partial_z f(z)={f(qz)-f(q^{-1}z) \over (q-q^{-1})z}$.
We have \cite{KS}

\begin{Theorem}: 
The Drinfeld generators at level-one are realized by the free boson fields as
\begin{eqnarray*}
& & \gamma=q,\\
& & h^1_m = a^1_mq^{-|m|/2}-a^{2}_mq^{|m|/2},
\ \ h^{2}_m =a^{2}_mq^{-|m|/2}+b_mq^{-|m|/2},~~~m\in{\bf Z}, 
 \\ 
& &X^{\pm}_1(z) = \pm :e^{\pm h_1(z;\pm\frac{1}{2})}:e^{\pm i\pi a^1_0},
 \ \ X^{+}_2(z) = :e^{h_{2}(z;\frac{1}{2})}e^{c(z;0)}: e^{-i\pi a^1_0}, 
 \\ 
& &X^{-}_2(z)= :e^{-h_{2}(z;-\frac{1}{2})}[ \partial_z
e^{-c(z;0)}]:e^{i\pi a^1_0}.
\end{eqnarray*}
\end{Theorem} 

\subsection{Level-one vertex operators}
Let $V$ be  the 3-dimensional  evaluation  representation of
$\Uq$ and $E_{ij}$ be the $3\times 3$ matrix whose $(i,j)$-element is unity
and zero  otherwise. The grading of the basis vectors
$v_1$, $v_2$, $v_3$ of $V$ is chosen to be $[v_1]=[v_2]=1$, $[v_3]=0$.Then 
the 3-dimenional level-0 representation $V_z$ of $\Uq$ is given by 
\begin{eqnarray*}
& &e_1=E_{12},\ \ e_2=E_{23},\ \ e_0=-zE_{31},\ \ t_1=q^{E_{11}-E_{22}},
\ \ 
t_2=q^{E_{22}+E_{33}}\\
& &f_1=E_{21},\ \ f_2=E_{32},\ \ f_0=z^{-1}E_{13},\ \
t_0=q^{-E_{11}-E_{33}},\ \ d=z\frac{d}{dz}.
\end{eqnarray*} 
We define the dual modules $V_z^{*S}$ of $V_z$ by 
$\pi_{V^{*S}}(a)=\pi_{V}(S(a))^{st}$, $\forall a\in \Uq$, where 
$st$ is the supertransposition operation.

Throughout, we denote by $V(\lambda)$ an irreducible
highest weight $\Uq$-module with 
highest weight $\lambda$. Consider the following intertwiners of 
$\Uq$-modules:
\begin{eqnarray*}
& &\Phi_{\lambda}^{\mu V}(z) :
 V(\lambda) \longrightarrow V(\mu)\otimes V_{z} ,\ \ \ \ 
\Phi_{\lambda}^{\mu V^{*}}(z) :
 V(\lambda) \longrightarrow V(\mu)\otimes V_{z}^{*S} ,\\
& &\Psi_{\lambda}^{V \mu}(z) :
 V(\lambda) \longrightarrow V_{z}\otimes V(\mu),\ \ \ \
\Psi_{\lambda}^{V^* \mu}(z) :
 V(\lambda) \longrightarrow V_{z}^{*S}\otimes V(\mu).
\end{eqnarray*}
They are intertwiners in the sense that for any $x\in \Uq$, 
\begin{eqnarray}
\label{EX}
\Theta(z)\cdot x=\Delta(x)\cdot \Theta(z),\ \ \ \Theta(z)=
\Phi(z),\Phi^{*}(z),\Psi(z),\Psi^{*}(z).
\end{eqnarray}
The intertwiners are even operators, that is their grading
is $[\Theta(z)]=0$. $\Phi(z)$ 
($\Phi^{*}(z)$) is called type I (dual) vertex operator and $\Psi(z)$ 
($\Psi^{*}(z)$) type II (dual) vertex operator.

We expand the vertex operators as 
\begin{eqnarray}
& &\Phi(z)=\sum_{j=1,2,3}\Phi_j(z)\otimes v_j\  ,\ \ \ \
\Phi^{*}(z)=\sum_{j=1,2,3}\Phi^{*}_j(z)\otimes v^{*}_j,\\
& &\Psi(z)=\sum_{j=1,2,3}v_j\otimes\Psi_j(z)\  ,\ \ \ \
\Psi^{*}(z)=\sum_{j=1,2,3}v^{*}_j\otimes\Psi^{*}_j(z).
\end{eqnarray}
Define the operators $\phi_j(z),\phi^{*}_j(z),\psi_j(z)$ and 
$\psi^{*}_j(z)$ $(j=1,2,3)$ by  
\bea
& &\phi_3(z)=:e^{-h^{*}_2(q^2z;-\frac{1}{2})+c(q^2z;0)}:
  e^{-i\pi  a^{2}_{0}},\no\\
& & \phi_2(z)=-[\phi_3(z),f_2]_{q^{-1}},\ \ \ \ 
  \phi_1(z)=[\phi_2(z),f_1]_q,\no\\
& &\phi^{*}_1(z)=:e^{h^{*}_1(qz;-\frac{1}{2})}:e^{i\pi a^{2}_{0}},\no\\ 
& &\phi^{*}_2(z)=-q^{-1}[\phi^{*}_1(z),f_1]_{q},\ \ \ \
  \phi^{*}_3(z)=q^{-1}[\phi^{*}_2(z),f_2]_q, \no\\ 
& &\psi_1(z)=:e^{-h^{*}_1(qz;\frac{1}{2})}:e^{-i\pi a^{2}_{0}},\no\\
& &\psi_2(z)=[\psi_1(z),e_1]_{q},\ \ \ \
  \psi_3(z)=[\psi_2(z),e_2]_q,\no\\
& &\psi^{*}_3(z)=:e^{h^{*}_2(z;\frac{1}{2})}_1
  \partial_z[e^{-c(z;0)}]:e^{i\pi a^{2}_{0}},\no\\
& &\psi^{*}_2(z)=q[\psi^{*}_3(z),e_2]_{q^{-1}},\ \ \ \
  \psi^{*}_1(z)=-q[\psi^{*}_2(z),e_1]_q,
\eea
where $h^{*1}_{m}=-h^2_m$,
$h^{*2}_m=-h^1_m-\frac{[2m]}{[m]}h^2_m$ and $Q_{h^{*1}}=-Q_{h^2}$, 
$ Q_{h^{*2}}=-Q_{h^1}-2Q_{h^2}$. Introduce 
$\phi(z),\phi^{*}(z),\psi(z),\psi^{*}(z)$ by
\bea
&&\phi(z)=\sum_{j=1,2,3}\phi_j(z)\otimes v_j\  ,\ \ \ \
  \phi^{*}(z)=\sum_{j=1,2,3}\phi^{*}_j(z)\otimes v^{*}_j,\label{VT1}\\
&&\psi(z)=\sum_{j=1,2,3}v_j\otimes\psi(z)_j\  ,\ \ \ \                 
  \psi^{*}(z)=\sum_{j=1,2,3}v^{*}_j\otimes\psi^{*}(z)_j.\label{VT2}   
\eea
Then we have 

\begin{Proposition} (\cite{KS}): The operators $\phi(z),\;\phi^*(z),\;
\psi(z)$ and $\psi^*(z)$ satisfy the same commutation relations as
$\Phi(z),\; \Phi^*(z),\; \Psi(z)$ and $\Psi^*(z)$, respectively.
\end{Proposition}

\section{Level-one highest weight $\Uq$-modules revisited}
In this section, we re-examine the level-one irreducible
highest weight $\Uq$-modules in the Fock
space defined by the bosonic $q$-oscillators.

Denote by $F_{\lambda_{a^1},\lambda_{a^2},\lambda_b,\lambda_c}$
the bosonic Fock spaces
generated by $a^{i}_{-m},b_{-m},c_{-m}$ $(m>0)$ over the 
vector $|\lambda_{a^1},\lambda_{a^2},\lambda_b,\lambda_c>$:
\begin{eqnarray}
 F_{\lambda_{a^1},\lambda_{a^2},\lambda_b,\lambda_c}={\bf
C}[a^{i}_{-1},a^{i}_{-2},...
;b_{-1},..;c_{-1},..]|\lambda_{a^1},\lambda_{a^2},\lambda_b,\lambda_c>,
\end{eqnarray}
 where 
\beq
|\lambda_{a^1},\lambda_{a^2},\lambda_b,\lambda_c>=e^{\lambda_{a^1}Q_{a^1}+
\lambda_{a^2}Q_{a^2}+\lambda_bQ_{b}+\lambda_cQ_c}|0>.
\eeq
The vacuum vector $|0>$ is defined by $a^i_m|0>=b_m|0>=c_m|0>=0$ for
$m>0$. Obviously, 
\begin{eqnarray*}
a^{i}_m|\lambda_{a^1},\lambda_{a^2},\lambda_b,\lambda_c>=0,\ \ 
{\rm for}\ \ m>0,\\
b_m|\lambda_{a^1},\lambda_{a^2},\lambda_b,\lambda_c>=0,\ \ {\rm for}
\ \ m>0,\\
c_m|\lambda_{a^1},\lambda_{a^2},\lambda_b,\lambda_c>=0, \ \ {\rm for }\ \ m>0.
\end{eqnarray*}
$|\l_{a^1},\l_{a^2},\l_b,\l_c>$ is said to be a $\Uq$ highest weight vector of
the weight $\l=\l^0\L_0+\l^1 \L_1+\l^2\L_2$, if it satisfies
\bea
e_i|\l_{a^1},\l_{a^2},\l_b,\l_c>&=&0,\no\\
h_i|\l_{a^1},\l_{a^2},\l_b,\l_c>&=&
   \l^i|\l_{a^1},\l_{a^2},\l_b,\l_c>,~~~i=0,1,2.
\eea
Here $\L_i~(i=0,1,2)$ are the
fundamental weights of $\Uq$.
In order to classify the highest weight $\Uq$-modules in the Fock
space, we introduce the spaces
\begin{eqnarray}
F_{(\alpha;\beta)}=\bigoplus_{i,j\in {\bf Z}} F
_{\beta+i,\beta-i+j,\beta-\alpha+j,-\alpha+j},
\end{eqnarray}
where $\a$ and $\b$ are arbitrary parameters. In the following we
restrict ourselves to $\a\in{\bf Z}$ case.

\noindent {\it Remark.} In \cite{KS}, the following two extra spaces
\bea
&&F_{((1,0);\beta)}
  =\bigoplus_{i,j\in {\bf Z}} F
  _{\beta+1+i,\beta-i+j,\beta+j,j},\no\\
&&F_{((0,1);\beta)} 
  =\bigoplus_{i,j\in {\bf Z}} F
  _{\beta+1+i,\beta+1-i+j,\beta+j,j}
\eea
were also introduced. However, it is easily seen that 
$F_{((1,0);\beta)},\; 
 F_{((0,1);\beta)}\subset F_{(\a;\beta)}$. In fact,
$F_{((1,0);\beta)}=F_{(1;\beta)}$ and 
 $F_{((0,1);\beta)}=F_{(2;\beta)}$.

It can be shown that the bosonized action of $\Uq$ on $F_{(\a;\b)}$
is closed, i.e. $\Uq F_{(\a;\b)}=F_{(\a;\b)}$. Hence each Fock space
$F_{(\alpha;\beta)}$ constitutes a $\Uq$-module.
However these modules are not irreducible in general.
To obtain the irreducible subspaces in $F_{(\alpha;\beta)}$, we
introduce a pair of fermionic currents\cite{KS}\cite{BCMN}
\begin{eqnarray*}
\eta(z)=\sum_{n\in{\bf Z}}\eta_nz^{-n-1}=:e^{c(z;0)}:,\ \ 
\xi(z)=\sum_{n\in{\bf Z}}\xi_nz^{-n}=:e^{-c(z;0)}:.
\end{eqnarray*}
The mode expansion of $\eta(z)$ and $\xi(z)$ is well defined on $
F_{(\alpha;\beta)}$ for $\alpha\in{\bf Z}$, and the modes satisfy the 
relations
\begin{eqnarray*}
\xi_m\xi_n+\xi_n\xi_m=\eta_m\eta_n+\eta_n\eta_m=0,\ \ 
\xi_m\eta_n+\eta_n\xi_m=\delta_{m,n}.
\end{eqnarray*}
Obviously, $\eta_0\xi_0$ and $\xi_0\eta_0$ qualify as projectors and so we
use them to decompose $F_{(\a;\b)}$ into a direct sum
$F_{(\a;\b)}=\eta_0\xi_0 F_{(\a;\b)}\oplus\xi_0\eta_0F_{(\a;\b)}$.
Following \cite{KS},
$\eta_0\xi_0F_{(\a;\b)}$ is referred to as Ker$_{\eta_0}$ and
$\xi_0\eta_0F_{(\a;\b)}=F_{(\a;\b)}/\eta_0\xi_0F_{(\a;\b)}$ as
Coker$_{\eta_0}$. 
Since $\eta_0$ commutes (or anticommutes) with $\Uq$, 
Ker$_{\eta_0}$ and Coker$_{\eta_0}$
are both $\Uq$-modules.

\subsection{Characters and supercharacters}
Now we study the characters and supercharacters of these
$\Uq$-modules Ker$_{\eta_0}$ and
Coker$_{\eta_0}$. We first of all bosonize
the derivation operator $d$ as
\begin{eqnarray}
\label{D}
& &d=\sum_{m=1}^\infty\frac{m^2}{[m]^2}\{h_{-m}^{1}h_{m}^{2}+h_{-m}^{2}h_{m}^{1}
+(q^m+q^{-m})h_{-m}^{2}h_{m}^{2}\}-\sum_{m=1}^\infty
\frac{m^2}{[m]^2}c_{-m}c_m,\nonumber\\
& &\ \ \ \ +\{h^{1}_0h^2_0+h^2_0h^2_0-\frac{1}{2}c_0(c_0+1)\}.
\end{eqnarray}
\noindent One can easily check that this $d$ obeys the
commutation relations,
\begin{eqnarray}
\label{GR}
[d,K^i]=0,\ \ [d,h^i_m]=mh^i_m,\ \ [d,X^{\pm,j}_m]=mX^{\pm,j}_m, 
\end{eqnarray}
as required. Moreover, we have $[d,\xi_0]=[d,\eta_0]=0$.

\noindent{\it Remark}. In\cite{KS}, the authors also gave a bosonized
expression for the operator $d$. However, their $d$ does not satisfy
the derivation properties (\ref{GR}).

The character and supercharacter of a $\Uq$-module $V$ are defined by 
\bea
ch_{V}(q,x,y)&=&tr_{V}(
q^{-2d}x^{h^1_0}y^{h^2_0}),\\
Sch_{V}(q,x,y)&=&Str_{V}(q^{-2d}x^{h^1_0}y^{h^2_0})
=tr_{V}((-1)^{N_f}q^{-2d}x^{h^1_0}y^{h^2_0}),
\eea
respectively. 
The Fermi-number operator $(-1)^{N_f}$ can also bosonized and we derive 
$(-1)^{N_f}=(-1)^{b_0}$. Indeed such a bosonized 
operator satisfies
\begin{eqnarray*}
(-1)^{b_0}\Theta(z)=(-1)^{[\Theta(z)]}\Theta(z)(-1)^{b_0},\ \ 
{\rm for} \ \
\Theta(z)=X^{\pm}_i(z),\Phi(z),\Phi^{*}(z),\Psi(z),\Psi^{*}(z),
\end{eqnarray*}
as required. The characters of Ker$_{\eta_0}$ and 
Coker$_{\eta_0}$ were claculated in \cite{KS} by inserting the 
projectors $\eta_0\xi_0$ and $\xi_0\eta_0$ into the trace over the Fock 
space $F_{(\a;\b)}$ and then computing the resultant traces. 
Here we rederive these character formulae by using the BRST
resolution. We also derive the supercharacters of
Ker$_{\eta_0}$ and 
Coker$_{\eta_0}$. 

Let us define the Fock spaces, for $l\in {\bf Z}$, 
\begin{eqnarray*}
F^{(l)}_{(\alpha;\beta)}
=\bigoplus_{i,j\in {\bf Z}} F
_{\beta+i,\beta-i+j,\beta-\alpha+j,-\alpha+j+l}.
\end{eqnarray*}
\noindent We have  $F^{(0)}_{(\alpha;\beta)}=F_{(\alpha;\beta)}$. 
It can be shown that $\eta_0,\xi_0$ intertwine various Fock spaces
\begin{eqnarray*}
& &\eta_0\ \ :\ \ F^{(l)}_{(\alpha;\beta)}\longrightarrow
F^{(l+1)}_{(\alpha;\beta)},\\
& &\xi_0\ \ :\ \ F^{(l)}_{(\alpha;\beta)}\longrightarrow
F^{(l-1)}_{(\alpha;\beta)}.
\end{eqnarray*}
\noindent We have the following BRST complexes:
\begin{eqnarray}
\label{GRA}
  \begin{array}{ccccccc}
    \cdot\cdot\cdot&\stackrel{Q_{l-1}=\eta_0}{\longrightarrow}&
    F^{(l)}_{(\alpha;\beta)}&\stackrel{Q_{l}=\eta_0}
    {\longrightarrow}&F^{(l+1)}_{(\alpha;\beta)}&\stackrel 
    {Q_{l+1}=\eta_0}{\longrightarrow}&\cdot\cdot\cdot\\  
    & &|{\bf O}&&|{\bf O}&&\\
    \cdot\cdot\cdot&\stackrel{Q_{l-1}=\eta_0}{\longrightarrow}&
    F^{(l)}_{(\alpha;\beta)}&\stackrel{Q_{l}=\eta_0}
    {\longrightarrow}&F^{(l+1)}_{(\alpha;\beta)}&\stackrel
    {Q_{l+1}=\eta_0}{\longrightarrow}&\cdot\cdot\cdot
  \end{array},
\end{eqnarray}
where ${\bf O}$ is an operator such that
$F^{(l)}_{(\alpha;\beta)}\longrightarrow F^{(l)}_{(\alpha;\beta)}$,
and ${\bf O}$ commutes with the BRST charges $Q_l$.
We have 

\begin{Proposition}: 
\begin{eqnarray}
& &Ker_{Q_{l}}=Im_{Q_{l-1}},\ \ {\rm for \ \ any \ \ }l\in{\bf Z} \ \ 
{\rm
and}\nonumber\\
& &tr({\bf O})|_{Ker_{Q_l}}=tr({\bf O})|_{Im_{Q_{l-1}}}=tr({\bf O})|_{
Coker_{Q_{l-1}}}.
\end{eqnarray}
\end{Proposition}
\noindent {\it Proof.} It follows from the fact that
$\eta_0\xi_0+\xi_0\eta_0=1$, $\xi_0^2=\eta_0^2=0$ and $\eta_0\xi_0$
($\xi_0\eta_0$ ) is the projection operator from $F_{(\a;\b)}$
to  $Ker_{\eta_0}$ 
$(Coker_{\eta_0})$.

In the following we simply write Ker$_{\eta_0}$ and Coker$_{\eta_0}$ of
$F_{(\a;\b)}$ as Ker$F_{(\a;\b)}$ and Coker$F_{(\a;\b)}$, respectively.
Then 
\begin{Proposition}: The characters of
$Ker_{F_{(\alpha;\beta)}}$ 
and $Coker_{F_{(\alpha;\beta)}}$ for $\alpha\in{\bf Z}$ are given by
\bea
ch_{Ker_{F_{(\alpha;\beta)}}}&=&ch_{Ker_{F^{(0)}_{(\alpha;\beta)}}}
=\sum_{l=1}^\infty(-1)^{l+1}ch_{F^{(-l)}_{(\alpha;\beta)}}=
\frac{q^{-\alpha(\alpha+1)}}{\prod_{n=1}^\infty(1-q^{2n})^3}\nonumber\\
& &\times \sum_{i,j\in Z}\sum_{l=1}^\infty(-1)^{l+1}
q^{l(l-1)+2l(\alpha-j)+(2i^2-2ij+j^2+j)}x^{2i-j}y^{\alpha-i},\no\\
ch_{Coker_{F_{(\alpha;\beta)}}}&=&ch_{Coker_{F^{(0)}_{(\alpha;\beta)}}}
=\sum_{l=1}^\infty(-1)^{l+1}ch_{F^{(l)}_{(\alpha;\beta)}}=
\frac{q^{-\alpha(\alpha+1)}}{\prod_{n=1}^\infty(1-q^{2n})^3}\nonumber\\
& &\times \sum_{i,j\in Z}\sum_{l=1}^\infty(-1)^{l+1}
q^{l(l+1)-2l(\alpha-j)+(2i^2-2ij+j^2+j)}x^{2i-j}y^{\alpha-i},
\eea
respectively.
\end{Proposition}

{\noindent \em Remark.} By using the following identity:
\begin{eqnarray*}
\sum_{l\in Z}(-1)^{l+1}q^{l(l-1)+2lt}=0,~~{\rm for~any}~~t\in{\bf Z}.
\end{eqnarray*}
one can show that the above character formulae coincide with those
given by  Kimura et al in \cite{KS}.

\begin{Proposition}: The supercharacters of
$Ker_{F_{(\alpha;\beta)}}$ and 
$CoKer_{F_{(\alpha;\beta)}}$ for $\alpha\in{\bf Z}$ are given by
\bea
Sch_{Ker_{F_{(\alpha;\beta)}}}&=&Sch_{Ker_{F^{(0)}_{(\alpha;\beta)}}}
=\sum_{l=1}^\infty(-1)^{l+1}Sch_{F^{(-l)}_{(\alpha;\beta)}}=
\frac{(-1)^{\beta-\alpha}q^{-\alpha(\alpha+1)}}
{\prod_{n=1}^\infty(1-q^{2n})^3}\nonumber\\
& &\times \sum_{i,j\in Z}(-1)^{j}\sum_{l=1}^\infty(-1)^{l+1}
q^{l(l-1)+2l(\alpha-j)+(2i^2-2ij+j^2+j)}x^{2i-j}y^{\alpha-i},\no\\
Sch_{Coker_{F_{(\alpha;\beta)}}}&=&Sch_{Coker_{F^{(0)}_{(\alpha;\beta)}}}
=\sum_{l=1}^\infty(-1)^{l+1}Sch_{F^{(l)}_{(\alpha;\beta)}}=
\frac{(-1)^{\beta-\alpha}q^{-\alpha(\alpha+1)}}
{\prod_{n=1}^\infty(1-q^{2n})^3}\nonumber\\
& &\times \sum_{i,j\in Z}(-1)^{j}\sum_{l=1}^\infty(-1)^{l+1}
q^{l(l+1)-2l(\alpha-j)+(2i^2-2ij+j^2+j)}x^{2i-j}y^{\alpha-i}.
\eea
\end{Proposition}

\noindent {\it Proof.} We sketch the proof of these two propositions.
Because $q^{-2d}x^{h^1_0}y^{h^2_0}$ and
$(-1)^{N_f}q^{-2d}x^{h^1_0}y^{h^2_0}$ commute with the BRST charges 
$Q_l$, the trace  over Ker and Coker can be written as 
the sum of trace  over $F^{(l)}_{(\alpha;\beta)}$. The latter
can be computed by the technique introduced in
\cite{CS}, which is given in appendix C.

Since $F_{(\alpha;\beta)}=F^{(1)}_{(\alpha-1;\beta+1)}$,
we have 

\begin{Corollary}: The following relations hold for any $\a$ and $\b$,
\begin{eqnarray*}
& &ch_{Coker_{F_{(\alpha;\beta)}}}= ch_{ker_{F_{(\alpha+1;\beta-1)}}} 
=ch_{ker_{F_{(\alpha+1;\beta)}}},\\
& &Sch_{Coker_{F_{(\alpha;\beta)}}}=Sch_{ker_{F_{(\alpha+1;\beta-1)}}}.
\end{eqnarray*}
\end{Corollary}

\subsection{$\Uq$-module structure of $F_{(\alpha;\beta-\alpha)}$}
Set $\l_\a=(1-\a)\L_0+\a\L_2$, and
\begin{eqnarray*} 
& &|\lambda_{\alpha}>=            
|\beta-\alpha,\beta-\alpha,\beta-2\alpha,-\alpha>~~\in
F_{(\alpha;\beta-\alpha)},~~~\a\in{\bf Z}, \\ 
& &|\Lambda_1>=|\beta,\beta-1,\beta-1,0>~~\in F_{(1;\beta-1)},\\
 & &|\Lambda_2>=|\beta-1,\beta-1,\beta-2,0>~~\in F_{(2;\beta-2)}.
\end{eqnarray*}
\noindent The above vectors play the role of  the highest weight vectors
of $\Uq$-modules \cite{KS}. 
One can check 
\bea
&&\eta_0|\lambda_{\alpha}>=0,~~ {\rm for}~\alpha=0,-1,-2,\cdots,\no\\
&&\eta_{0}|\Lambda_1>=0,~~~~\eta_0|\L_2>=0,\no\\
&&\eta_0|\lambda_{\alpha}>\ne 0,~~ {\rm for}~ \alpha=1,2,\cdots.\label{el} 
\eea
It follows that the modules Coker$F_{(\alpha;\beta-\alpha)}$ 
($\alpha=1,2,3,...$), Ker$F_{(\alpha;\beta-\alpha)}$ 
($\alpha=0,-1,-2,...$), Ker$F(1;\b-1)$ and Ker$F(2,\b-2)=
{\rm Coker}F(1;\b-1)$ are highest weight $\Uq$-modules.
Denote by $\bar{V}(\l_\a)$ and $\bar{V}(\L_1)$ the $\Uq$-modules with
highest weights $\l_\a$ and $\L_1$, respectively.
From (\ref{el}) and corollary 1, we find 

\begin{Theorem}\label{v=ker}: We have the following identifications
of the highest weight $\Uq$-modules:
\bea
\bar{V}(\l_\a)&\cong&
  {\rm Ker}F_{(\a;\b-\a)}\equiv{\rm Coker}F_{(\a-1;\b-\a+1)},~~~
  {\rm for}~~\alpha=0,-1,-2,...,\no\\
&\cong&{\rm Coker}F_{(\a;\b-\a)}\equiv{\rm Ker}F_{(\a+1;\b-\a-1)},~~~
  {\rm for}~~\a=1,2,\cdots,\no\\
\bar{V}(\Lambda_1)&\cong& {\rm Ker}F_{(1;\b-1)}\equiv{\rm Coker}F_{(0;\b)}.
\eea
\end{Theorem}
Therefore we have

\begin{Theorem}: For $\a\in{\bf Z}$, each Fock space
$F_{(\alpha;\beta-\alpha)}$ can 
be decomposed explicitly into a direct sum of the highest weight $\Uq$-modules
\begin{eqnarray*}
 \begin{array}{cccc}
  &Ker && Coker\\
  .&.   && . \\
  .& .&&.\\
F_{(-2;\beta+2)}=&\bar{V}(\lambda_{-2})&\bigoplus& \bar{V}(\lambda_{-1})\\
&&&\\&\phi(z)\uparrow\downarrow \phi^{*}(z)&
     &\phi(z)\uparrow\downarrow \phi^{*}(z)\\&&&\\
F_{(-1;\beta+1)}=&\bar{V}(\lambda_{-1})&\bigoplus& \bar{V}(\Lambda_0)\\
 &&&\\ &\phi(z)\uparrow\downarrow \phi^{*}(z)&   
     &\phi(z)\uparrow\downarrow \phi^{*}(z)\\&&&\\
F_{(0;\beta)}=&\bar{V}(\Lambda_0)&\bigoplus& \bar{V}(\Lambda_1)\\
&&&\\     &\phi(z)\uparrow\downarrow \phi^{*}(z)&   
     &\phi(z)\uparrow\downarrow \phi^{*}(z)\\&&&\\
F_{(1;\beta-1)}=&\bar{V}(\Lambda_1)&\bigoplus& \bar{V}(\Lambda_2)\\
 &&&\\    &\phi(z)\uparrow\downarrow \phi^{*}(z)&   
     &\phi(z)\uparrow\downarrow \phi^{*}(z)\\&&&\\
F_{(2;\beta-2)}=&\bar{V}(\Lambda_2)&\bigoplus& \bar{V}(\lambda_2)\\
  &&&\\   &\phi(z)\uparrow\downarrow \phi^{*}(z)&   
     &\phi(z)\uparrow\downarrow \phi^{*}(z)\\&&&\\
F_{(3;\beta-3)}=&\bar{V}(\lambda_2)&\bigoplus& \bar{V}(\lambda_3)\\
  .& .&&.\\
  .& .&&.
 \end{array}
 \ \ .
\end{eqnarray*}
\end{Theorem}

It is expected that the modules $\bar{V}(\l_\a)$ and $\bar{V}(\L_1)$ are also 
irreducible with respect to the action
of $\Uq$. From theorem \ref{v=ker}. we are led to the following
conjecture which is a modified version of the conjecture
2.2 proposed in \cite{KS}.
\vskip.1in
\noindent{\bf Conjecture 1}: {\it $\bar{V}(\l_\a)$ and $\bar{V}(\L_1)$
are irreducible highest weight $\Uq$-modules of the weights $\l_\a$
and $\L_1$, respectively, i.e.
\bea
&&\bar{V}(\l_\a)=V(\l_\a),\ \ \ \ \a\in{\bf Z},\no\\
&&\bar{V}(\L_1)=V(\L_1).
\eea
    }

\section{Exchange relations of vertex operators}
In this section, we derive the exchange relations of the type I and type II 
vertex operators of $\Uq$. As expected, these vertex operators 
satisfy the graded Faddeev-Zamolodchikov algebra.

\subsection{The R-matrix}
Let $R(z)\in End(V\otimes V)$
be the R-matrix of $\Uq$: 
\begin{eqnarray}
\label{R}
R(z)(v_i\otimes v_j)=\sum_{kl}R^{ij}_{kl}(z)v_k\otimes v_l, \ \ \ \
 \forall v_i, v_j, v_k, v_l\in V,
\end{eqnarray}
where the matrix elements are given by
\begin{eqnarray*}
& &R^{33}_{33}(\frac{z_1}{z_2})=-\frac{z_1q^{-1}-z_2q}{z_1q-z_2q^{-1}},\ \ 
R^{23}_{23}(\frac{z_1}{z_2})=-\frac{z_1-z_2}{z_1q-z_2q^{-1}},\ \
R^{32}_{23}(\frac{z_1}{z_2})=\frac{(q-q^{-1})z_2}{z_1q-z_2q^{-1}},\\
& &R^{32}_{32}(\frac{z_1}{z_2})=-\frac{z_1-z_2}{z_1q-z_2q^{-1}},\ \
R^{23}_{32}(\frac{z_1}{z_2})=\frac{(q-q^{-1})z_1}{z_1q-z_2q^{-1}},\ \
R^{22}_{22}(\frac{z_1}{z_2})=-1,\\ 
& &R^{13}_{13}(\frac{z_1}{z_2})=-\frac{z_1-z_2}{z_1q-z_2q^{-1}},\ \
R^{31}_{13}(\frac{z_1}{z_2})=\frac{(q-q^{-1})z_2}{z_1q-z_2q^{-1}},\ \
R^{31}_{31}(\frac{z_1}{z_2})=-\frac{z_1-z_2}{z_1q-z_2q^{-1}},\\ 
& &R^{13}_{31}(\frac{z_1}{z_2})=\frac{(q-q^{-1})z_1}{z_1q-z_2q^{-1}},\ \
R^{12}_{12}(\frac{z_1}{z_2})=-\frac{z_1-z_2}{z_1q-z_2q^{-1}},\ \
R^{21}_{12}(\frac{z_1}{z_2})=-\frac{(q-q^{-1})z_2}{z_1q-z_2q^{-1}},\\ 
& &R^{21}_{21}(\frac{z_1}{z_2})=-\frac{z_1-z_2}{z_1q-z_2q^{-1}},\ \
R^{12}_{21}(\frac{z_1}{z_2})=-\frac{(q-q^{-1})z_1}{z_1q-z_2q^{-1}},\ \
R^{11}_{11}(\frac{z_1}{z_2})=-1 ,\\
& &R^{ij}_{kl}=0\ \ ,\ \ {\rm otherwise}.
\end{eqnarray*} 
The R-matrix satisfies  the graded Yang-Baxter equation(YBE) on 
$V\otimes V\otimes V$
\begin{eqnarray*}
R_{12}(z)R_{13}(zw)R_{23}(w)=R_{23}(w)R_{13}(zw)R_{12}(z),
\end{eqnarray*}
and moreover enjoys: (i) initial condition,
$R(1)=P$ with $P$ being
the graded permutation operator; (ii) unitarity condition,
$R_{12}(\frac{z}{w})R_{21}(\frac{w}{z})=1$, where 
$R_{21}(z)=PR_{12}(z)P$; and (iii) crossing-unitarity, 
\begin{eqnarray*}
R^{-1,st_1}(z)\lt((q^{-2\overline{\rho}}\otimes 1)R(zq^{-2})
(q^{2\overline{\rho}}\otimes 1)\rt)^{st_1}=1\otimes 1 ,
\end{eqnarray*}
\noindent where 
\begin{eqnarray}
\label{CU}
 q^{2\overline{\rho}}\stackrel{def}{=}
\left(
   \begin{array}{ccc}
     q^{2\rho_1}&&\\&q^{2\rho_2}&\\&&q^{2\rho_3}
    \end{array}
  \right)
=\left(
    \begin{array}{ccc}
     1&&\\&q^{-2}&\\&&q^{-2}
    \end{array}
  \right) .
\end{eqnarray}
The various supertranspositions of the R-matrix are given by
\begin{eqnarray*}
& &(R^{st_1}(z))^{kl}_{ij}=R(z)^{il}_{kj}(-1)^{[i]([i]+[k])},\ \ 
(R^{st_2}(z))^{kl}_{ij}=R(z)^{kj}_{il}(-1)^{[j]([j]+[l])},\\
& &(R^{st_{12}}(z))^{kl}_{ij}=
R(z)^{ij}_{kl}(-1)^{([i]+[j])([i]+[j]+[l]+[k])}
=R(z)^{ij}_{kl} .
\end{eqnarray*}

\subsection{The graded Faddeev-Zamolodchikov algebra}
Now, we are in the position to calculate the exchange relations of the type I 
and type II vertex operators of $\Uq$.
Define 
\begin{eqnarray*}
\oint dzf(z)=Res(f)=f_{-1} \ \ ,\ \ {\rm for \ \ formal\ \  series \ \ 
function\ \ } f(z)=\sum_{n\in Z}f_nz^{n} .
\end{eqnarray*}
\noindent Then the Chevalley generators of $\Uq$ can be expressed 
by the integrals,
\begin{eqnarray*}
& &e_1=\oint dzX^+_1(z),\ \ e_2=\oint dzX^{+}_2(z),\ \ 
  f_1=\oint dz X^-_1(z),\ \ f_2=\oint dzX^-_2(z), \\
& &e_0=-\oint\oint dzdw\;z\,[X^-_2(w),X^-_1(z)]_{q^{-1}}q^{-h^1_0-h^2_0},\\
& &f_0=\oint\oint
  dzdw\;z^{-1}q^{h^1_0+h^2_0}[X^+_1(z),X^+_2(w)]_{q}.
\end{eqnarray*}
One can also obtain the integral expression of the vertex 
operators defined in (\ref{VT1})-(\ref{VT2})
\bea
\phi_2(z)&=&\oint dw :\{\frac{e^{-c(wq;0)}}{wq(1-\frac{qz}{w})}
 +\frac{e^{-c(wq^{-1};0)}}{zq^2(1-\frac{w}{zq^3})}\}
 e^{-h^{*}_2(q^2z;-\frac{1}{2})-h_2(w;-\frac{1}{2})+c(q^2z;0)}
 e^{i\pi h^1_0}: ,\no\\
\phi_1(z)&=&\oint dw_1\oint dw
\frac{q^2-1}{w(1-\frac{w_1q}{w})(1-\frac{wq}{w_1})}\no\\
& & \times :\{\frac{e^{-c(wq;0)}}{wq(1-\frac{qz}{w})}
  +\frac{e^{-c(wq^{-1};0)}}{zq^2(1-\frac{w}{zq^3})}\}
 e^{-h^{*}_2(q^2z;-\frac{1}{2})-h_2(w;-\frac{1}{2})
 -h_1(w_1;-\frac{1}{2})+c(q^2z;0)}
 e^{-i\pi a^2_0}: ,\no\\
\phi^{*}_2(z)&=&\oint dw\frac{1-q^{-2}}{z(1-\frac{zq^2}{w})
 (1-\frac{w}{z})}:e^{h^{*}_1(qz;-\frac{1}{2})-h_1(w;-\frac{1}{2})}
 e^{-i\pi h^{1}_0}: ,\no\\
\phi^{*}_3(z)&=&\oint dw_1\oint dw \frac{1-q^{-2}}{z(1-\frac{zq^2}{w})
 (1-\frac{w}{z})}\no\\
& & \times :\frac{e^{-c(w_1q;0)}-e^{-c(w_1q^{-1};0)}}
 {ww_1(1-\frac{wq}{w_1})(1-\frac{w_1q}{w})}
 e^{h^{*}_1(qz;-\frac{1}{2})-h_1(w;-\frac{1}{2})-h_2(w_1;-\frac{1}{2})}
 e^{i\pi a^{2}_0}: ,\no\\
\psi_2(z)&=&\oint dw\frac{1-q^{2}}{wq(1-\frac{w}{zq^2})
 (1-\frac{z}{w})}:e^{-h^{*}_1(qz;\frac{1}{2})+h_1(w;\frac{1}{2})}
 e^{i\pi h^{1}_0}: ,\no\\
\psi_3(z)&=&\oint dw_1\oint dw
 \frac{(q^2-1)(q-q^{-1})}{wq(1-\frac{w}{zq^2})
 (1-\frac{z}{w})}\no\\
& & \times \frac{
 :e^{-h^{*}_1(qz;\frac{1}{2})+h_1(w;\frac{1}{2})+h_2(w_1;\frac{1}{2})
 +c(w;0)}e^{-i\pi a^{2}_0}:}
 {w_1(1-\frac{w_1}{wq})(1-\frac{w}{w_1q})},\no\\
\psi^{*}_2(z)&=&\oint dw\{ 
 \frac{e^{-c(zq;0)}}{w(1-\frac{zq}{w})}+
 \frac{e^{-c(zq^{-1};0)}}{zq^{-1}(1-\frac{wq}{z})}\}
 :e^{h^{*}_2(z;\frac{1}{2})+h_2(w;\frac{1}{2})+c(w;0)}
 e^{-i\pi h^{1}_0}: ,\no\\
\psi^{*}_1(z)&=&\oint dw\oint dw_1\frac{1-q^2}
 {w_1(1-\frac{w}{w_1q})(1-\frac{w_1}{wq})}\no\\
& &\times :\{
 \frac{e^{-c(zq;0)}}{w(1-\frac{zq}{w})}+
 \frac{e^{-c(zq^{-1};0)}}{zq^{-1}(1-\frac{wq}{z})}\}
 e^{h^{*}_2(z;\frac{1}{2})+h_2(w;\frac{1}{2})+h_1(w_1;\frac{1}{2})+c(w;0)}
 e^{i\pi a^{2}_0}: .
\eea

Using these integral expressions and
the relations given in appendix A and appendix B, we derive

\begin{Proposition}: The bosonic vertex operators defined in 
(\ref{VT1})--(\ref{VT2}) satisfy the graded Faddeev-Zamolodchikov
algebra
\begin{eqnarray}
& &\phi_j(z_2)\phi_i(z_1)=\sum_{kl}R(\frac{z_1}{z_2})^{kl}_{ij}
\phi_k(z_1)\phi_l(z_2)(-1)^{[i][j]} ,\\
& &\phi^{*}_j(z_2)\phi^{*}_i(z_1)=\sum_{kl}R(\frac{z_1}{z_2})^{ij}_{kl}
\phi^{*}_k(z_1)\phi^{*}_l(z_2)(-1)^{[i][j]},\\
& &\psi_i(z_1)\psi_j(z_2)=\sum_{kl}R(\frac{z_1}{z_2})^{kl}_{ij}
\psi_l(z_2)\psi_k(z_1)(-1)^{[i][j]} ,\\
& &\psi^{*}_i(z_1)\psi^{*}_j(z_2)=\sum_{kl}R(\frac{z_1}{z_2})^{ij}_{kl}
\psi^{*}_l(z_2)\psi^{*}_k(z_1)(-1)^{[i][j]} ,\\
& &\psi^{*}_i(z_1)\phi_j(z_2)=-(-1)^{[i][j]}\phi_j(z_2)\psi^{*}_i(z_1),\\
& &\phi_j(z_2)\phi^{*}_i(z_1)=\sum_{kl}\overline{R}(\frac{z_1}{z_2})^{kl}_{ij}
  \phi^{*}_k(z_1)\phi_l(z_2)(-1)^{[k][l]},
\end{eqnarray}
where $\bar{R}(z)=R^{-1,st_1}(z)$. 
\end{Proposition}

In the derivation of this proposition the fact that
$R(z)^{kl}_{ij}(-1)^{[k][l]}$ = $R(z)^{kl}_{ij}(-1)^{[i][j]}$ is
helpful.

\begin{Proposition}: We have the first invertibility relations,
\bea
&&\phi_i(z)\phi^{*}_j(z)=-(-1)^{[j]}\delta_{ij},\no\\
&&-\sum_{k}(-1)^{[k]}\phi^{*}_k(z)\phi_k(z)=1,
\eea
and the second invertibility relations,
\bea
&&\phi^{*}_i(zq^2)\phi_j(z)=\delta_{ij}q^{2\rho_i} ,\no\\
&&\sum_{k}q^{-2\rho_k}\phi_k(z)\phi^{*}_k(zq^2)=1 .
\eea
\end{Proposition}

Using the fact that $\eta_0\xi_0$ is a projection operator, 
we can make the following identifications:
\bea
&&\Phi_i(z)=\eta_0\xi_0\phi_i(z)\eta_0\xi_0,~~~~
\Phi^{*}_i(z)=\eta_0\xi_0\phi^{*}_i(z)\eta_0\xi_0 ,\label{V1}\\
&&\Psi_i(z)=\eta_0\xi_0\psi_i(z)\eta_0\xi_0,~~~~
\Psi^{*}_i(z)=\eta_0\xi_0\psi^{*}_i(z)\eta_0\xi_0  .\label{V2}
\eea
Since the vertex operators $\phi,\phi^{*},\psi,\psi^{*}$
commute (or anti-commute) with the BRST charge $\eta_0$, we have 

\begin{Theorem}: Set
\beq
\mu_\a=\lt\{
\begin{array}{ll}
\L_\a, & ~~\a=0,1,2,\\
\l_{\a-1}, & ~~{\rm for}~\a>2,\\
\l_\a, & ~~{\rm for}~\a<0.
\end{array}
\rt.
\eeq
Then the vertex operators defined by
(\ref{V1}) and (\ref{V2}) interwine the level-one irreducible
highest weight $\Uq$-modules
$V(\mu_{\alpha})~~(\alpha\in{\bf Z})$
\begin{eqnarray*}
& &\Phi(z) :
 V(\mu_{\alpha}) \longrightarrow V(\mu_{\alpha-1})
\otimes V_{z} ,\ \
\Phi^{*}(z) :
 V(\mu_{\alpha}) \longrightarrow V(\mu_{\alpha+1})\otimes
V_{z}^{*S} ,\\
& &\Psi(z) :
 V(\mu_{\alpha}) \longrightarrow V_{z}\otimes V(\mu_{\alpha-1})
,\ \
\Psi^{*}(z) :
 V(\mu_{\alpha}) \longrightarrow V_{z}^{*S}\otimes
V(\mu_{\alpha+1}),   
\end{eqnarray*}
and satisfy the graded Faddeev-Zamolodchikov algebra.
Moreover, the type I vertex operators $\Phi(z),\Phi^{*}(z)$ obey the
following invertibility relations:
\begin{eqnarray}
& &\Phi_i(z)\Phi^{*}_j(z)|_{V(\mu_{\alpha})}
=-(-1)^{[j]}\delta_{ij}id_{V(\mu_{\alpha})} ,\\
& &-\sum_{k}(-1)^{[k]}\Phi^{*}_k(z)\Phi_k(z)
|_{V(\mu_{\alpha})}=id|_{V(\mu_{\alpha)}} ,\nonumber\\
& &\Phi^{*}_i(zq^2)\Phi_j(z)|_{V(\mu_{\alpha})}
=\delta_{ij}q^{2\rho_i}id|_{V(\mu_{\alpha})},\\
& &\sum_{k}q^{-2\rho_k}\Phi_k(z)\Phi^{*}_k(zq^2)
|_{V(\mu_{\alpha})}=id|_{V(\mu_{\alpha})} .\nonumber 
\end{eqnarray}  
\end{Theorem}

\subsection{Miki's construction of $\Uq$}

We generalize the Miki's construction to the supersymmetric case.
Define 
\begin{eqnarray*}
& &(L^+(z))^{j}_{i}=\Phi_i(zq^{\frac{1}{2}})
\Psi^{*}_j(zq^{-\frac{1}{2}}), \\ 
& &(L^-(z))^{j}_{i}=\Phi_i(zq^{-\frac{1}{2}})\Psi^{*}_j(zq^{\frac{1}{2}}).
\end{eqnarray*}

\begin{Proposition}:  The L-operators $L^{\pm}(z)$ defined above give a
realization of the super Reshetikhin-Semenov-Tian-Shansky algebra at 
level one in the quantum affine superalgebra $\Uq$ \cite{Z,GZ}
\begin{eqnarray*}
& &R(\frac{z}{w})L^{\pm}_1(z)L^{\pm}_2(w)=
L^{\pm}_2(w)L^{\pm}_1(z)R(\frac{z}{w}),\\
& &R(\frac{z^+}{w^-})L^{+}_1(z)L^{-}_2(w)=
L^{-}_2(w)L^{+}_1(z)R(\frac{z^-}{w^+}) ,
\end{eqnarray*}
where $L^{\pm}_1(z)=L^{\pm}(z)\otimes 1$, 
$L^{\pm}_2(z)=1\otimes L^{\pm}(z)$ and $z^{\pm}=zq^{\pm\frac{1}{2}}$.
\end{Proposition}

\noindent{\em Proof.} Straightforward by means of the graded
Faddeev-Zamolodchikov algebra.

\section{$q$-deformed supersymmetric $t$-$J$ model }
In this section, we give a mathematical definition of the $q$-deformed
supersymmetric $t$-$J$ model on an infinite lattice.

\subsection{Space of states}
By means of the R-matrix (\ref{R}) of $\Uq$,
one can define the $q$-deformed
supersymmetric $t$-$J$ model on the infinte lattice $\cdots
\otimes V\otimes V\otimes V\cdots$. Let $h$ be the operator on
$V\otimes V$ such that
\begin{eqnarray*}
& &PR(\frac{z_1}{z_2})=1+u h+\cdots,\ \ \ \ \ \ \ \ u
\longrightarrow 0 ,\\
& &\ \ \ \ \  P:{\rm the\ \  graded\ \ permutation \ \ operator },\ \
 e^{u}\equiv \frac{z_1}{z_2}.
\end{eqnarray*}
The Hamiltonian $H$ of the $q$-deformed supersymmetric $t$-$J$ model
is defined by
\begin{eqnarray}
H=\sum_{l\in Z}h_{l+1,l}.
\end{eqnarray}
$H$ acts formally on the infinite tensor product,
\beq
\cdots V\otimes V\otimes V\cdots.\label{vvv}
\eeq
It can be easily checked that 
\begin{eqnarray*}
[U'_q(\widehat{sl}(2|1)), H]=0,
\end{eqnarray*}
where $U'_q(\widehat{sl}(2|1))$ is the subalgebra of $\Uq$ with the
derivation operator $d$ being dropped. So $U'_q(\widehat{sl}(2|1))$ plays
the role of infinite dimensional {\it non-abelian symmetries} of the
$q$-deformed supersymmetric $t$-$J$ model on the infinite lattice. 
Following \cite{DFJMN}, we replace the infinite tensor product
(\ref{vvv}) by the level-0 $\Uq$-module,
\begin{eqnarray*}
{\it F}_{\a\a'}={\rm Hom}(V(\mu_\a),V(\mu_{\a'}))\cong V(\mu_\a)\otimes
  V(\mu_{\a'})^{*},
\end{eqnarray*}
where $V(\mu_\a)$ is level-one irreducible highest weight
$\Uq$-module and $V(\mu_{\a'})^{*}$ is the dual module of $V(\mu_{\a'})$.
By theorem 4, this homomorphism can be realized by applying the type I vertex
operators repeatedly.
So we shall make the (hypothetical) identification:
\begin{eqnarray*}
``{\bf \rm the \ \ space \ \ of \ \ physical \ \ states } " =
\bigoplus_{\a,\a'\in Z} V (\mu_\a)\otimes V(\mu_{\a'})^{*}.
\end{eqnarray*}
Namely, we take
\begin{eqnarray*}
F\equiv End(\bigoplus_{\a\in Z}V(\mu_\a))\cong 
 \bigoplus_{\a,\a'\in{\bf Z}} F_{\a\a'}
\end{eqnarray*}
as the space of states of the $q$-deformed supersymmetric $t$-$J$
model on the infinite lattice.
The left action of $\Uq$ on ${\it F}$ is defined by
\begin{eqnarray*}
x.f=\sum x_{(1)}\circ f\circ S(x_{(2)})(-1)^{[f][x_{(2)}]},~~~
 \forall x\in\Uq,~f\in F,
\end{eqnarray*}
where we have used notation $\Delta(x)=\sum x_{(1)}
\otimes x_{(2)}$.
Note that
$F_{\a\a}$ has the unique canonical element $id_{V(\mu_\a)}$. We call
it the vacuum \cite{JM} and denote it by
$|vac>_{\a}\in {\it F}_{\a\a}$.

\subsection{Local structure and local operators}
Following Jimbo et al\cite{JM}, we use the type I 
vertex operators and their variants to incorporate the local structure into the
space of physical states ${\it F}$, that is to formulate the action
of local operators of the $q$-deformed supersymmetric $t$-$J$ model on
the infinite tensor product (\ref{vvv}) in terms of their actions on
$F_{\a\a'}$. 

Using the isomorphisms (c.f. theorem 4)
\bea
\Phi(1)&:&~V(\mu_\a)\longrightarrow V(\mu_{\a-1})\otimes V,\no\\
\Phi^{*,st}(q^2)&:&~V\otimes V(\mu_\a)^*\longrightarrow
  V(\mu_{\a-1})^*,
\eea
were $st$ is the supertransposition on the quantum space,
we have the following identification:
\begin{eqnarray*}
V(\mu_\a)\otimes V(\mu_{\a'})^{*}\stackrel{\sim}{\rightarrow}
V(\mu_{\a-1})\otimes V\otimes V(\mu_{\a'})^{*}\stackrel{\sim}{\rightarrow}
V(\mu_{\a-1})\otimes V(\mu_{\a'-1})^{*}.
\end{eqnarray*}
The resulting isomorphism can be identified with
the super translation (or shift) operator defined by
\begin{eqnarray*}
T=\sum_i\Phi_i(1)\otimes \Phi_i^{*,st}(q^2)(-1)^{[i]}q^{-2\rho_i}.
\end{eqnarray*}
Its inverse is given by
\begin{eqnarray*}
T^{-1}=\sum_i\Phi_i^{*}(1)\otimes \Phi_i^{st}(1).
\end{eqnarray*}

Thus we can define the local operators on $V$ as operators on
$F_{\a\a'}$ \cite{JM}. Let us label the tensor components from the middle as
$1,2,\cdots$ for the left half and as $0,-1,-2,\cdots$ for
the right half. The operators acting on the site 1  are defined by
\begin{eqnarray}
E_{ij}\stackrel{def}{=}E^{(1)}_{ij}=
-\Phi^{*}_i(1)\Phi_j(1)(-1)^{[j]}\otimes id.
\end{eqnarray}
More generally we set
\begin{eqnarray}
\label{L}
E^{(n)}_{ij}=T^{-(n-1)}E_{ij}T^{n-1}\ \ \ \ (n\in Z).
\end{eqnarray}
Then, from the invertibility relations of the type I vertex operators
of $\Uq$, we have

\begin{Theorem}: The local operators $E^{(n)}_{ij}$ acting on
${\it F}_{\a\a'}$  satisfy the following relations:
\begin{eqnarray*}
E^{(m)}_{ij}E^{(n)}_{kl}=\left\{
  \begin{array}{ll}
  \delta_{jk}E^{(n)}_{il} & {\rm if}\ \ \ m=n\\
  (-1)^{([i]+[j])([k]+[l])}E^{(n)}_{kl}E^{(m)}_{ij}&{\rm if}\ \ \ m\ne n
  \end{array}
 \right. .
\end{eqnarray*}
\end{Theorem}

As is expected from the physical point of view, we also have

\begin{Proposition}: The vacuum vectors
$|vac>_{\a}$ are super-translationally invariant and singlets (i.e. belong
to the trivial representation of $\Uq$)
\begin{eqnarray*}
T|vac>_\a=|vac>_{\a-1} ,\\
x.|vac>_\a=\epsilon(x)|vac>_\a .
\end{eqnarray*}
\end{Proposition}
\noindent {\it Proof.} Let $u^{(\a)}_l$ ($u^{*(\a)}_l$) be a basis vectors of
$V(\mu_\a)\ \  (V(\mu_\a)^{*})$  and
\begin{eqnarray*}
|vac>_\a\stackrel{def}{=}id_{V(\mu_\a)}=\sum_l u^{(\a)}_l\otimes u^{*(\a)}_l.
\end{eqnarray*}
Then
\beq
T|vac>_\a=\sum_{m,l}q^{-2\rho_m}\Phi_m(1)u^{(\a)}_l\otimes \Phi^{*,st}_m
(q^2)u^{*(\a)}_l(-1)^{[m]+[l][m]}.
\eeq
We want to show $T|vac>_\a=|vac>_{\a-1}$. This
is equivalent to proving
\begin{eqnarray*}
\sum_{m,l}q^{-2\rho_m}\Phi_m(1)u^{(\a)}_l\,\Phi^{*,st}_m(q^2)\cdot
  u^{*(\a)}_l(v)(-1)
  ^{[m]+[l][m]}=v,\ \ {\rm for \ \ any \ \ } v\in V(\mu_{\a-1}).
\end{eqnarray*}
Now
\bea
 l.h.s&=&\sum_{m,l}q^{-2\rho_m}\Phi_m(1)u^{(\a)}_l u^{*(\a)}_l\lt
(\Phi^{*}_m(q^2)^{st})^{st}v\rt)(-1)^{[m]}\no\\
&=&\sum_{m,l}q^{-2\rho_m}\Phi_m(1)u^{(\a)}_l u^{*(\a)}_l
  (\Phi^{*}_m(q^{2})v)\no\\
&=&\sum_{m}q^{-2\rho_m}\Phi_m(1)\Phi^{*}_m(q^2)v=v ,\no
\eea
where we have used $(\Phi^{*}_m(z)^{st})^{st}
=\Phi^{*}_m(z)(-1)^{[m]}$ and the second invertibility of
the type I vertex operators. As to the second equation, we have
\bea
x\cdot|vac>_\a&=&\sum x_{(1)}u^{(\a)}_l\otimes x_{(2)}u^{*(\a)}_l
  (-1)^{[l][x_{(2)}]}\no\\
&=&\sum x_{(1)}u^{(\a)}_l\otimes \pi_{V(\mu_\a)^*}(x_{(2)})_{lm}u^{*(\a)}_m
  (-1)^{[l][x_{(2)}]}\no\\
&=&\sum x_{(1)}u^{(\a)}_l\otimes \pi_{V(\mu_\a)}
  (S(x_{(2)}))_{ml}u^{*(\a)}_m\no\\
&=&\sum x_{(1)}\pi_{V(\mu_\a)}(S(x_{(2)}))_{ml}u^{(\a)}_l\otimes u^{*(\a)}_m
  \no\\
&=&\sum x_{(1)}S(x_{(2)})u^{(\a)}_m\otimes u^{*(\a)}_m
  =\e(x)|vac>_\a.\no
\eea
This completes the proof.

For any local operator $O\in F$, its vacuum expectation value
is given by
\begin{eqnarray*}                                                           
{}_\a<vac|O|vac>_\a=
  \frac{tr_{V(\mu_\a)}(q^{-2d+2h^2_0}O)}
  {tr_{V(\mu_\a)}(q^{-2d+2h^2_0})},
\end{eqnarray*}
where we have chosen the normalization ${}_\a<vac|vac>_\a=1$. We shall
denote the correlator ${}_\a<vac|O|vac>_\a$ by $<O>_\a$.

\section{Correlation functions} 
The aim of  this section is to calculate $<E_{mn}>_{\a}$.
The generalization to the calculation of the  multi-point functions 
is straightforward.

Throughout we use the abbreviation
\begin{eqnarray}
(z;x)_\infty=\prod^\infty_{n=0}(1-zx^n).
\end{eqnarray}
Set 
\begin{eqnarray*}
P^m_n(z_1,z_2|q|\a)=
\frac{tr_{V(\mu_\a)}(q^{-2d+2h^2_0}
\Phi^{*}_m(z_1)\Phi_n(z_2))}
{tr_{V(\mu_\a)}(q^{-2d+2h^2_0})} ,
\end{eqnarray*}
then $<E_{mn}>_{\a}=P^m_n(z,z|q|\a)$. By (\ref{V1})--(\ref{V2}),
conjecture 1 and theorem 2, it is sufficient to calculate
\begin{eqnarray}
F^{(\alpha)}_{mn}(z_1,z_2)=
\frac{tr_{F_{(\alpha;\beta-\alpha)}}(q^{-2d+2h^2_0}
\phi^{*}_m(z_1)\phi_n(z_2)\eta_0\xi_0)}
{tr_{F_{(\alpha;\beta-\alpha)}}(q^{-2d+2h^2_0}\eta_0\xi_0)}.
\end{eqnarray}
By the procedure similar to the derivation of the (super)characters, we get
\begin{eqnarray*}
F^{(\alpha)}_{mn}(z_1,z_2)=\frac{\delta_{mn}}{\chi_{\alpha}}
\sum_{l=1}^\infty(-1)^{l+1}F_{m,-l},
\end{eqnarray*}
where 
\bea
\chi_{\alpha}&=&tr_{F_{(\alpha;\beta-\alpha)}}(q^{-2d+2h^2_0}\eta_0\xi_0)
=\frac{q^{-\alpha(\alpha+1)}}{\prod_{n=1}^\infty(1-q^{2n})^3}\no\\
& &\times \sum_{i,j\in Z}\sum_{l=1}^\infty(-1)^{l+1}
  q^{l(l-1)+2l(\alpha-j)+(2i^2-2ij+j^2+j)+2(\alpha-i)},\no\\
F_{1,l}&=&
\frac{q^{-\alpha(\alpha+1)}(\frac{z_2}{z_1}q^2;q^2)_\infty
(\frac{z_1}{z_2}q^2;q^2)_\infty}{(q^2;q^2)_\infty}\sum_{i,j\in Z}
q^{l(l+1)+2l(j-\alpha)+2i^2-2ij+j^2+j}\no\\
& &\times \oint dw_1\oint dw 
\frac{1-q^2}{w^2q}\frac{(\frac{z_2}{w}q^5;q^2)_\infty(\frac{q^2z_2}{w_1})^{2i-j}
(\frac{z_2}{z_1}q^4)^{\alpha-i}(\frac{qz_2}{w})^{l-i+j}}
{(\frac{w_1}{w}q;q^2)_\infty(\frac{w}{w_1}q;q^2)_\infty
(\frac{w_1}{z_1};q^2)_\infty
(\frac{z_1}{w_1}q^4;q^2)_\infty(\frac{z_2}{w}q;q^2)_\infty}\no\\
& &\times \oint dw_1\oint dw
\frac{1-q^2}{wz_2q^2}\frac{(\frac{w}{z_2}q;q^2)_\infty
(\frac{q^2z_2}{w_1})^{2i-j}
(\frac{z_2}{z_1}q^{2})^{\alpha-i}(\frac{q^3z_2}{w})^{l-i+j}}  
{(\frac{w_1}{w}q;q^2)_\infty(\frac{w}{w_1}q;q^2)_\infty
(\frac{w_1}{z_1};q^2)_\infty
(\frac{z_1}{w_1}q^4;q^2)_\infty(\frac{w}{z_2}q^{-3};q^2)_\infty} ,\no\\
F_{2,l}&=&-
\frac{q^{-\alpha(\alpha+1)}(\frac{z_2}{z_1}q^2;q^2)_\infty
(\frac{z_1}{z_2}q^2;q^2)_\infty}{(q^2;q^2)_\infty}\sum_{i,j\in Z}
q^{l(l+1)+2l(j-\alpha)+2i^2-2ij+j^2+j}\no\\
& &\times \oint dw_1\oint dw
\frac{q-q^{-1}}{w_1wq}\frac{(\frac{z_2}{w}q^5;q^2)_\infty
(\frac{q^2z_2}{w_1})^{2i-j}
(\frac{z_2}{z_1}q^4)^{\alpha-i}(\frac{qz_2}{w})^{l-i+j}}
{(\frac{w_1}{w}q^3;q^2)_\infty(\frac{w}{w_1}q;q^2)_\infty
(\frac{w_1}{z_1};q^2)_\infty
(\frac{z_1}{w_1}q^2;q^2)_\infty(\frac{z_2}{w}q;q^2)_\infty}\no\\
& &\times \oint dw_1\oint dw
\frac{q-q^{-1}}{w_1z_2q^2}\frac{(\frac{w}{z_2}q;q^2)_\infty
(\frac{q^2z_2}{w_1})^{2i-j}
(\frac{z_2}{z_1}q^{2})^{\alpha-i}(\frac{q^3z_2}{w})^{l-i+j}}
{(\frac{w_1}{w}q^3;q^2)_\infty(\frac{w}{w_1}q;q^2)_\infty
(\frac{w_1}{z_1};q^2)_\infty
(\frac{z_1}{w_1}q^2;q^2)_\infty(\frac{w}{z_2}q^{-3};q^2)_\infty},\no\\
F_{3,l}&=&
\frac{q^{-\alpha(\alpha+1)}(\frac{z_2}{z_1}q^2;q^2)_\infty
(\frac{z_1}{z_2}q^2;q^2)_\infty}{(q^2;q^2)_\infty}\sum_{i,j\in Z}
q^{l(l+1)+2l(j-\alpha)+2i^2-2ij+j^2+j}\no\\
& &\times \oint dw_1\oint dw
\frac{q^{-1}-q}{w_1wq}\frac{(\frac{z_2}{w}q^3;q^2)_\infty
(\frac{q^2z_2}{w_1})^{2i-j}
(\frac{z_2}{z_1}q^4)^{\alpha-i}(\frac{qz_2}{w})^{l-i+j}}
{(\frac{w_1}{w}q;q^2)_\infty(\frac{w}{w_1}q;q^2)_\infty
(\frac{w_1}{z_1};q^2)_\infty
(\frac{z_1}{w_1}q^2;q^2)_\infty(\frac{z_2}{w}q;q^2)_\infty}\no\\
& &\times \oint dw_1\oint dw
\frac{q^2-1}{ww_1}\frac{(\frac{w}{z_2}q;q^2)_\infty
(\frac{q^2z_2}{w_1})^{2i-j}
(\frac{z_2}{z_1}q^{2})^{\alpha-i}(\frac{q^3z_2}{w})^{l-i+j}}
{(\frac{w_1}{w}q;q^2)_\infty(\frac{w}{w_1}q;q^2)_\infty
(\frac{w_1}{z_1};q^2)_\infty
(\frac{z_1}{w_1}q^2;q^2)_\infty(\frac{w}{z_2}q^{-1};q^2)_\infty}.
\eea

We now derive the difference equations satisfied by these one-point
functions. Let
\begin{eqnarray*}
\overline{F}^{(\alpha)}_{mn}(z_1,z_2)=
tr_{F_{(\alpha;\beta-\alpha)}}(q^{-2d+2h^2_0}
\phi^{*}_m(z_1)\phi_n(z_2)\eta_0\xi_0)\stackrel{def}{=}\delta_{mn}
\overline{F}^{(\alpha)}_m(z_1,z_2).
\end{eqnarray*}
Noticing that 
\begin{eqnarray*}
& &x^{d}\phi_i(z)x^{-d}=\phi_i(zx^{-1}),\ \ 
x^{d}\phi^{*}_i(z)x^{-d}=\phi^{*}_i(zx^{-1}),\\
& &x^{d}\psi_i(z)x^{-d}=\psi_i(zx^{-1}),\ \
x^{d}\psi^{*}_i(z)x^{-d}=\psi^{*}_i(zx^{-1}),\\
& &x^{d}\eta_0x^{-d}=\eta_0,\ \ 
x^{d}\xi_0x^{-d}=\xi_0,
\end{eqnarray*}
we get the difference equations
\begin{eqnarray*}
\overline{F}^{(\alpha)}_m(z_1,z_2q^2)=q^{-2\rho_m}\,\sum_{k}R(z_2,z_1)^{km}_{mk}
 \overline{F}^{(\alpha-1)}_k(z_1,z_2)
 \,(-1)^{[m]+[k]+[m][k]}.
\end{eqnarray*}
Since $\a\in{\bf Z}$, it is easily seen that this is a set of infinite
number of difference equations.

\section*{Acknowledgements}
This work was financially supported by the Australian Research Council.
The first named author would like to thank the second named author
and department of mathematics, the
University of Queensland, for their kind hospitality. The first
named author was also  partially
supported by the National Natural Science Foundation of China and Northwest
University Fund.

\section*{Appendix A.}
In this apppendix, we give the normal order relations of fundmental bosonic
fields:
\begin{eqnarray*}
& &e^{h_1(z_1;\beta_1)}e^{h_1(z_2;\beta_2)}=(z_1-q^{-(\beta_1+\beta_2)+1}z_2)
(z_1-q^{-(\beta_1+\beta_2)-1}z_2):e^{h_1(z_1;\beta_1)}e^{h_1(z_2;\beta_2)}:,\\
& &e^{h_1(z_1;\beta_1)}e^{h_2(z_2;\beta_2)}=\frac{1}
{z_1-q^{-(\beta_1+\beta_2)}z_2}:e^{h_1(z_1;\beta_1)}e^{h_2(z_2;\beta_2)}:,\\
& &e^{h_2(z_1;\beta_1)}e^{h_1(z_2;\beta_2)}=\frac{1}
{z_1-q^{-(\beta_1+\beta_2)}z_2}:e^{h_2(z_1;\beta_1)}e^{h_1(z_2;\beta_2)}:,\\
& &e^{h_2(z_1;\beta_1)}e^{h_2(z_2;\beta_2)}=
:e^{h_2(z_1;\beta_1)}e^{h_2(z_2;\beta_2)}:,\\
& &e^{h_i(z_1;\beta_1)}e^{h^{*}_j(z_2;\beta_2)}=
(z_1-q^{-(\beta_1+\beta_2)}z_2)^{\delta_{ij}}:
e^{h_i(z_1;\beta_1)}e^{h^{*}_j(z_2;\beta_2)}:,\\
& &e^{h^{*}_i(z_1;\beta_1)}e^{h_j(z_2;\beta_2)}=
(z_1-q^{-(\beta_1+\beta_2)}z_2)^{\delta_{ij}}:
e^{h^{*}_i(z_1;\beta_1)}e^{h_j(z_2;\beta_2)}:,\\
& &e^{h^{*}_1(z_1;\beta_1)}e^{h^{*}_1(z_2;\beta_2)}=
:e^{h^{*}_1(z_1;\beta_1)}e^{h^{*}_1(z_2;\beta_2)}:,\\
& &e^{h^{*}_1(z_1;\beta_1)}e^{h^{*}_2(z_2;\beta_2)}=\frac{1}
{z_1-q^{-(\beta_1+\beta_2)}z_2}:e^{h^{*}_1(z_1;\beta_1)}
e^{h^{*}_2(z_2;\beta_2)}:,\\
& &e^{h^{*}_2(z_1;\beta_1)}e^{h^{*}_1(z_2;\beta_2)}=\frac{1}
{z_1-q^{-(\beta_1+\beta_2)}z_2}:e^{h^{*}_2(z_1;\beta_1)}
e^{h^{*}_1(z_2;\beta_2)}:,\\
& &e^{h^{*}_2(z_1;\beta_1)}e^{h^{*}_2(z_2;\beta_2)}=\frac{1}
{(z_1-q^{-(\beta_1+\beta_2)+1}z_2)(z_1-q^{-(\beta_1+\beta_2)-1}z_2)}
:e^{h^{*}_2(z_1;\beta_1)}e^{h^{*}_2(z_2;\beta_2)}:,\\
& &e^{c(z_1;\beta_1)}e^{c(z_2;\beta_2)}=(z_1-q^{-(\beta_1+\beta_2)}z_2)
:e^{c(z_1;\beta_1)}e^{c(z_2;\beta_2)}:.
\end{eqnarray*}
\section*{Appendix B. }
By means of the bosonic realiztion of $\Uq$, the integral expressions of
the vertex operators and the technique given in Ref.\cite{AJ}, one 
can check the following relations

\begin{itemize}
\item For the type I vertex operators
\begin{eqnarray*}
& &[\phi_1(z),f_1]_{q^{-1}}=0,\ \ \phi_1(z)=[\phi_2(z),f_1]_{q},\ \ 
[\phi_3(z),f_1]=0,\\
& &[\phi_1(z),f_2]=0,\ \ [\phi_2(z),f_2]_{q^{-1}}=0,\ \
\phi_2(z)=-[\phi_3(z),f_2]_{q^{-1}},\\           
& &[\phi_1(z),e_1]=t_1\phi_2(z),\ \ [\phi_2(z),e_1]=0,\ \
[\phi_3(z),e_1]=0,\\           
& &[\phi_1(z),e_2]=0,\ \ -[\phi_2(z),e_2]=t_2\phi_3(z),\ \
[\phi_3(z),e_2]=0,\\
& &\phi_1(z)t_1=qt_1\phi_1(z),\ \ \phi_2(z)t_1=q^{-1}t_1\phi_2(z),\ \ 
\phi_3(z)t_1=t_1\phi_3(z),\\
& &\phi_1(z)t_2=t_2\phi_1(z),\ \ \phi_2(z)t_2=qt_2\phi_2(z),\ \ 
\phi_3(z)t_2=qt_2\phi_3(z),
\end{eqnarray*}
\begin{eqnarray*}
& &\phi^{*}_2(z)=-q^{-1}[\phi^{*}_1(z),f_1]_{q},\ \
[\phi^{*}_2(z),f_1]_{q^{-1}}=0,\ \
[\phi^{*}_3(z),f_1]=0,\\
& &[\phi^{*}_1(z),f_2]=0,\ \
\phi^{*}_3(z)=q^{-1}[\phi^{*}_2(z),f_2]_{q},\ \
[\phi^{*}_3(z),f_2]_{q}=0,\\
& &[\phi^{*}_1(z),e_1]=0,\ \ [\phi^{*}_2(z),e_1]=-q^{-1}
t_1\phi_1^{*}(z),\ \ [\phi^{*}_3(z),e_1]=0,\\   
& &[\phi^{*}_1(z),e_2]=0,\ \ -[\phi^{*}_2(z),e_2]=0,\ \
[\phi^{*}_3(z),e_2]=-q^{-1}t_2\phi^{*}_2(z),\\
& &\phi^{*}_1(z)t_1=q^{-1}t_1\phi^{*}_1(z),\ \
\phi^{*}_2(z)t_1=qt_1\phi^{*}_2(z),\ \
\phi^{*}_3(z)t_1=t_1\phi^{*}_3(z),\\
& &\phi^{*}_1(z)t_2=t_2\phi^{*}_1(z),\ \
\phi^{*}_2(z)t_2=q^{-1}t_2\phi^{*}_2(z),\ \
\phi^{*}_3(z)t_2=q^{-1}t_2\phi^{*}_3(z).
\end{eqnarray*}
\item For the type II vertex operators
\begin{eqnarray*}
& &\psi_2(z)=[\psi_1(z),e_1]_{q},\ \ [\psi_2(z),e_1]_{q^{-1}}=0,\ \
[\psi_3(z),e_1]=0,\\
& &[\psi_1(z),e_2]=0,\ \ \psi_3(z)=[\psi_2(z),e_2]_{q},\ \
[\psi_3(z),e_2]_{q}=0,\\
& &[\psi_1(z),f_1]=0,\ \ [\psi_2(z),f_1]=t^{-1}_1\psi_1(z),\ \
[\psi_3(z),f_1]=0,\\   
& &[\psi_1(z),f_2]=0,\ \ [\psi_2(z),f_2]=0,\ \
[\psi_3(z),f_2]=t^{-1}_2\psi_2(z),\\
& &\psi_1(z)t_1=qt_1\psi_1(z),\ \ \psi_2(z)t_1=q^{-1}t_1\psi_2(z),\ \
\psi_3(z)t_1=t_1\psi_3(z),\\
& &\psi_1(z)t_2=t_2\psi_1(z),\ \ \psi_2(z)t_2=qt_2\psi_2(z),\ \
\psi_3(z)t_2=qt_2\psi_3(z),
\end{eqnarray*}
\begin{eqnarray*}
& &[\psi^{*}_1(z),e_1]_{q^{-1}}=0,\ \
\psi^{*}_1(z)=-q[\psi^{*}_2(z),e_1]_{q},\ \
[\psi^{*}_3(z),e_1]=0,\\
& &[\psi^{*}_1(z),e_2]=0,\ \ [\psi^{*}_2(z),e_2]_{q^{-1}},\ \
\psi^{*}_2(z)=q[\psi^{*}_3(z),e_2]_{q^{-1}},\\
& &[\psi^{*}_1(z),f_1]=-qt^{-1}_1\psi^{*}_2(z),\ \
[\psi^{*}_2(z),f_1]=0,\ \
[\psi^{*}_3(z),f_1]=0,\\
& &[\psi^{*}_1(z),f_2]=0,\ \ [\psi^{*}_2(z),f_2]=-qt^{-1}_2\psi^{*}_3(z)
,\ \ [\psi^{*}_3(z),f_2]=0,\\
& &\psi^{*}_1(z)t_1=q^{-1}t_1\psi^{*}_1(z),\ \
\psi^{*}_2(z)t_1=qt_1\psi^{*}_2(z),\ \
\psi^{*}_3(z)t_1=t_1\psi^{*}_3(z),\\
& &\psi^{*}_1(z)t_2=t_2\psi^{*}_1(z),\ \
\psi^{*}_2(z)t_2=q^{-1}t_2\psi^{*}_2(z),\ \
\psi^{*}_3(z)t_2=q^{-1}t_2\psi^{*}_3(z).
\end{eqnarray*}
\end{itemize}

\section*{Appendix C.}
In computating the characters/supercharacters and the correlation 
functions, one encounters the trace of the form
\begin{eqnarray*}
tr(x^{-d}e^{\sum_{m=1}^\infty(\sum_{i=1,2}A_{m,i}h^i_{-m})+B_mc_{-m}}
e^{\sum_{m=1}^\infty(\sum_{i=1,2}C_{m,i}h^i_{m})+D_mc_{m}}f_1^{a^1_0}
f_2^{a^2_0}f_3^{b_0}f_4^{c_0}),   
\end{eqnarray*}
where $A_{m,i},B_m,C_{m,i},D_m$ and $f_i$ are all some
coefficients. 
We can calculate the contributions 
from the oscillators modes and the zero modes {\it separately}. The trace over
the oscillator modes can be carried out as follows 
by using the Clavelli-Shapiro 
technique \cite{CS}.
Let us introduce the extra oscillators 
$h'^i_m ,c'_m$ which commutate with $h^i_m,c_m$. 
$h'^i_m ,c'_m$ satisfy the same commutation relations as those 
satisfied by $h^i_m,c_m$. 
Introduce the operators 
\begin{eqnarray*}
& &H^{i}_m=\frac{h^i_m\otimes 1}{1-x^m}+1\otimes h'^i_m,\ \ \ \ 
C_m=\frac{c_m\otimes 1}{1-x^m}+1\otimes c'_m  ,\ \ m>0 ,\\
& &H^{i}_m=h^i_m\otimes 1+\frac{1\otimes h'^i_m}{x^m-1},\ \ \ \ 
C_m=c_m\otimes 1+\frac{1\otimes c'_m}{x^m-1} ,\ \ m<0. 
\end{eqnarray*}
Then for any bosonic operator $O(h^i_m,c_m)$, one can show 
\begin{eqnarray*}
tr(x^{-d}O(h^i_m,c_m))=\frac{
<0|O(H^i_m,C_m)|0>
}{\prod _{n=1}^\infty(1-x^n)^3} 
\end{eqnarray*}
providing that $d$ satisfies the derivation properties (\ref{GR}).
We write $<0|O(H^i_m,C_m)|0>\equiv <<O(h^i_m,c_m)>>$. Then by the
Wick theorem, one obtains
\begin{eqnarray*}
& &<<:e^{h_1(z_1;-\frac{1}{2})}e^{h_2(z_2;-\frac{1}{2})}:>>
=C_1(-\frac{1}{2})C_2(-\frac{1}{2})g_1(\frac{z_2}{z_1})\ \ ,\\
& &<<:e^{h_1(z_1;-\frac{1}{2})}e^{h_2(z_2;-\frac{1}{2})}:>>
=C_1(-\frac{1}{2})C_2(-\frac{1}{2})g_2(\frac{z_2}{z_1})\ \ ,\\
& &<<:e^{h_2(z_1;-\frac{1}{2})}e^{h_1(z_2;-\frac{1}{2})}:>>
=C_1(-\frac{1}{2})C_2(-\frac{1}{2})g_2(\frac{z_2}{z_1})\ \ ,\\
& &<<:e^{h_2(z_1;-\frac{1}{2})}e^{h_2(z_2;-\frac{1}{2})}:>>
=C_2(-\frac{1}{2})C_2(-\frac{1}{2})\ \ ,\\
& &<<:e^{h_i(z_1;-\frac{1}{2})}e^{h^{*}_j(z_2;-\frac{1}{2})}:>>
=C_i(-\frac{1}{2})C^{*}_j(-\frac{1}{2})
(g_2(\frac{z_2}{z_1}))^{-\delta_{ij}}\ \ ,\\
& &<<:e^{h^{*}_i(z_1;-\frac{1}{2})}e^{h_j(z_2;-\frac{1}{2})}:>>
=C^{*}_i(-\frac{1}{2})C_j(-\frac{1}{2})
(g_2(\frac{z_2}{z_1}))^{-\delta_{ij}}\ \ ,\\
& &<<:e^{h^{*}_1(z_1;-\frac{1}{2})}e^{h^{*}_1(z_2;-\frac{1}{2})}:>>
=C^{*}_1(-\frac{1}{2})C^{*}_1(-\frac{1}{2})\ \ ,\\
& &<<:e^{h^{*}_1(z_1;-\frac{1}{2})}e^{h^{*}_2(z_2;-\frac{1}{2})}:>>
=C^{*}_1(-\frac{1}{2})C^{*}_2(-\frac{1}{2})g_2(\frac{z_2}{z_1})\ \ ,\\
& &<<:e^{h^{*}_2(z_1;-\frac{1}{2})}e^{h^{*}_1(z_2;-\frac{1}{2})}:>>
=C^{*}_1(-\frac{1}{2})C^{*}_2(-\frac{1}{2})g_2(\frac{z_2}{z_1})\ \ ,\\
& &<<:e^{h^{*}_2(z_1;-\frac{1}{2})}e^{h^{*}_2(z_2;-\frac{1}{2})}:>>
=C^{*}_2(-\frac{1}{2})C^{*}_2(-\frac{1}{2})
(g_1(\frac{z_2}{z_1}))^{-1}\ \ ,\\
& &<<:e^{c(z_1;0)}e^{c(z_2;0)}:>>
=C_0(0)C_0(0)g_0(\frac{z_2}{z_1})\ \ ,
\end{eqnarray*}
\noindent where 
\begin{eqnarray*}
& &g_1(z)=(xzq^2;x)_\infty
  (xz;x)_\infty(xz^{-1}q^2;x)_\infty(xz^{-1};x)_\infty,\ \ \ 
  g_2(z)=\frac{1}{(xzq;x)_\infty(xz^{-1}q;x)_\infty},\\
& &g_0(z)=(xz;x)_\infty(xz^{-1};x)_\infty,\ \ \ 
   C_0(0)=(x;x)_\infty,\\
& & C_1(-1/2)=(xq^2;x)_\infty(x;x)_\infty=\frac{1}{C^{*}_2(-1/2)},\ \ \ 
 C_2(-1/2)=\frac{1}{C^{*}_1(-1/2)}=1.
\end{eqnarray*}

\end{document}